\documentclass[11pt,leqno]{amsart}
\usepackage{amssymb,verbatim,enumerate,ifthen}
\usepackage[mathscr]{eucal}
\usepackage[utf8]{inputenc}
\usepackage[T1]{fontenc}
\usepackage{cite}
\def\N{\mathbb{N}}
\def\R{\mathbb{R}}

\def\C{\mathscr{C}}

\newtheorem{theorem}{Theorem}[section]
\newtheorem*{theorem*}{Theorem}
\def\Thm#1#2{\ifthenelse{\equal{#1}{*}}{\begin{theorem*}#2\end{theorem*}}
             {\begin{theorem}\label{T#1}#2\end{theorem}}}
\newtheorem{Atheorem}{Theorem}

\def\thm#1{Theorem~\ref{T#1}}
\newtheorem{proposition}[theorem]{Proposition}
\newtheorem*{proposition*}{Proposition}
\def\Prp#1#2{\ifthenelse{\equal{#1}{*}}{\begin{proposition*}#2\end{proposition*}}
{\begin{proposition}\label{P#1}#2\end{proposition}}}

\newtheorem{corollary}[theorem]{Corollary}
\newtheorem*{corollary*}{Corollary}
\def\Cor#1#2{\ifthenelse{\equal{#1}{*}}{\begin{corollary*}#2\end{corollary*}}
             {\begin{corollary}\label{C#1}#2\end{corollary}}}

\newtheorem{lemma}[theorem]{Lemma}
\newtheorem*{lemma*}{Lemma}
\def\Lem#1#2{\ifthenelse{\equal{#1}{*}}{\begin{lemma*}#2\end{lemma*}}
             {\begin{lemma}\label{L#1}#2\end{lemma}}}
\def\lem#1{Lemma~\ref{L#1}}

\theoremstyle{definition}
\newtheorem{remark}[theorem]{Remark}
\newtheorem*{remark*}{Remark}
\def\Rem#1#2{\ifthenelse{\equal{#1}{*}}{\begin{remark}\rm #2\end{remark}}
             {\begin{remark}\label{R#1}\rm #2\end{remark}}}

\newtheorem{example}[theorem]{Example}
\newtheorem*{example*}{Example}
\def\Exa#1#2{\ifthenelse{\equal{#1}{*}}{\begin{example*}\rm #2\end{example*}}
             {\begin{example}\label{Ex#1}\rm #2\end{example}}}

\def\eq#1{{\rm(\ref{E#1})}}
\def\Eq#1#2{\ifthenelse{\equal{#1}{*}}
  {\begin{equation*}\begin{aligned}#2\end{aligned}\end{equation*}}
  {\begin{equation}\begin{aligned}\label{E#1}#2\end{aligned}\end{equation}}}

\begin{document}
\begin{flushright}
%\textit{Submitted to:} 
\end{flushright}
\vspace{5mm}

\date{\today}

\title{On approximately Convex and Affine Sequences}

\author[A. R. Goswami]{Angshuman R. Goswami}
\address[A. R. Goswami]{Department of Mathematics, University of Pannonia, 
H-8200 Veszprem, Hungary}
\email{goswami.angshuman.robin@mik.uni-pannon.hu}

\begin{abstract}
In this paper, our primary objective is to study a possible decomposition of an approximately convex sequence. \\

For a given $\varepsilon>0$; a sequence $\big<u_n\big>_{n=0}^{\infty}$ is said to be $\varepsilon$-convex, if for any $i,j\in\N$ with $i<j$ there exists an $n\in]i,j]\cap \N$ such that the following discrete functional inequality holds
\Eq{*}{
u_i-u_{i-1}-\dfrac{\varepsilon}{n-i}\leq u_j-u_{j-1}.
}
We show that such a sequence can be represented as the algebraic summation of a convex and a controlled sequence which is bounded in between $\left[-\dfrac{\varepsilon}{2}, \dfrac{\varepsilon}{2}\right].$\\

On the other hand, if for any $i,j\in\N$ with $i<j$, if a sequence $\big<u_n\big>_{n=0}^{\infty}$ satisfies the following form of inequality
\Eq{*}{
\left|\big(u_i-u_{i-1}\big)-\big(u_j-u_{j-1}\big)\right|\leq\dfrac{\varepsilon}{n-i}\quad \quad\mbox{for some}\quad n\in]i,j]\cap\N;
}
then we term it as $\varepsilon$-affine sequence. Such a sequence can be decomposed as the algebraic summation of an affine and a bounded sequence whose supremum norm doesn't exceed $\varepsilon.$\\

Various definitions, backgrounds, motivations, findings, and other important things are discussed in the introduction section.
\end{abstract}

\maketitle

\section*{Introduction} 
Throughout this paper $\N$, $\R$, and $\R_+$ will denote the set of natural, real, and non-negative numbers respectively.

Instead of treating sequences as ordered multi-sets;  we can visualize sequences as discrete functions. There are several well-known sequences that often appears in the study of various branches of mathematics. Arithmetic, Geometric, Fibonacci, Tribonacci, and Partitioning are some of the familiar examples. Most of these sequences can be studied under the umbrella of convex sequences. Most probably for the very first time, the terminology of sequential convexity was used in the book of
\cite{Mitrinovicc}. A sequence $\big<u_n\big>_{n=0}^{\infty}$ is termed as a convex sequence if the following inequality holds.
\Eq{99}{
2u_n\leq u_{n-1}+u_{n+1}\quad\mbox{for all} \quad n\in\N.
}  
In the case of reversed inequality, we will say $\big<u_n\big>_{n=0}^{\infty}$ is sequential concave. From the definition, it is visible that $\big<u_n\big>_{n=0}^{\infty}$ is a convex sequence if and only if the following inequality is satisfied
\Eq{*}{
u_i-u_{i-1}\leq u_j-u_{j-1}\quad\mbox{for all}\quad i,j\in\N \quad \mbox{with} \quad i<j.
}
In other words, monotonicity(increasingness) of the sequence $\big<u_n-u_{n-1}\big>_{n=1}^{\infty}$ implies convexity of $\big<u_n\big>_{n=0}^{\infty}$ and vice-verse. 
It is clear that if a sequence is both convex and concave; then it is nothing but an arithmetic sequence (arithmetic progression). The recent findings in the direction of sequential convexity can be found in the papers  
\cite{ Debnath,GauSte,Latreuch,pecaric,Essen,Jim} etc. and the references therein.\\

Convex sequences are also termed discrete versions of convex function. Let $I\subseteq\R$ be a non-empty interval. A function $f:I\to\R$ is said to be convex, if for any $x,y\in I$ and for all $t\in[0,1]$, the following inequality holds
\Eq{*}{
f(tx+(1-t)y)\leq tf(x)+(1-t) f(y).
}
Many researchers have studied various forms of convex functions such as higher-order convexity, approximate convexity, strong convexity etc. Many of the findings in functional convexity prove immense possibility to discover similar results in sequential convexity.  For instance, the discrete version of Hermite-Hadamard inequality.  Inspired by these various formulations of convexity, recently many mathematicians started to work in the direction of higher order convex sequences. Some of the interesting findings can be found in the papers  \cite{Krasniqi}, \cite{HDCS}. But, so far not many advancements have been made where a convex sequence is actually treated in an approximate sense. \\

Our primary objective in this paper is to define an approximate convex sequence and establish a stability theorem that resembles the famous Hyers-Ulam's stability theorem of convex functions.\\

In 1952, Hyers and Ulam in their paper \cite{Hyers} showed that for a given $\varepsilon>0$; if a function $f:I\to\R$ satisfies the following functional inequality
\Eq{*}{
f(tx+(1-t)y)\leq tf(x)+(1-t)f(y)+\varepsilon \quad\mbox{for any} \quad x,y\in I \quad \mbox{for all}\quad \mbox \quad t\in[0,1];
}
then there exists a convex function $C:I\to\R$ that closely approximate $f$ and satisfies \newline
$||f-C||\leq\varepsilon/2.$ We are going to use this result extensively in our paper.\\

A sequence $\big<u_n\big>_{n=0}^{\infty}$ is said to be $\varepsilon$-convex, if for any pair of $i,j\in\N$ with $i<j$; there exists atleast one $n\in [i+1,j]\cap\N$ such that the following inequality holds
\Eq{1}{
u_i-u_{i-1}\leq u_j-u_{j-1}+\dfrac{\varepsilon}{n-i}.
}
The set of sequences that satisfy \eq{1} is a subclass of a more generalized class of approximate convex sequences which contains elements $\big<u_n\big>_{n=0}^{\infty}$ that holds the following functional inequality 
\Eq{*}{
u_i-u_{i-1}\leq u_j-u_{j-1}+{\varepsilon}\quad \quad \mbox{where}\quad i,j\in \N\quad \mbox{with} \quad i<j.
}
In one of our recently submitted papers, we investigated about it.  It is clearly evident that any ordinary convex sequence also possesses $\varepsilon$-convexity. Besides three consecutive sequential points always form an approximate convex triplet. For verification,  we can substitute $j=i+1$ in $\eq{1}$ and we get the following form of discrete functional inequality 
\Eq{*}{
2u_i\leq u_{i-1}+u_{i+1}+{\varepsilon} \quad \mbox{for all}\quad i\in \N.
}
We showed that if a sequence $\big<u_n\big>_{n=0}^{\infty}$ is $\varepsilon$-convex, then we can express it as an algebraic summation of two very specific types of sequences  $\big<v_n\big>_{n=0}^{\infty}$ and $\big<w_n\big>_{n=0}^{\infty}$. The sequence $\big<v_n\big>_{n=0}^{\infty}$ possesses sequential convexity, while $\big<w_n\big>_{n=0}^{\infty}$ is a bounded sequence that fluctuate between $[-\varepsilon/2, \varepsilon/2].$\\

On the other hand, if both the $\big<u_n\big>_{n=0}^{\infty}$ and $\big<-u_n\big>_{n=0}^{\infty}$ are $\varepsilon$-convex, we term such a sequence as  $\varepsilon$-affine. We show that a $\varepsilon$-affine sequence can be approximated by an arithmetic sequence. \\

At the very end, we discussed an alternative definition of the convex sequence. Let $I\subseteq\R$ be a non-empty and non-singleton interval.  In the field of function theory, a function $f:I\to\R$ is said to be Wright convex if for any $x,y\in I$ and for all $t\in[0,1]$; the following functional inequality is satisfied
\Eq{*}{
f(tx+(1-t)y)+f((1-t)x+ty)\leq f(x)+f(y).
}
In other words, for any $[u,v]\subset[x,y]\subseteq I.$ with $u+v=x+y;$ a function $f:I\to\R$ is labeled as Wright convex if the following functional inequality holds
\Eq{*}{
f(u)+f(v)\leq f(x)+f(y).
}
The definition of Wright convexity was introduced  by E.M. Wright in 1954 in his publication \cite{Wright}.
It can be easily shown that any convex function $f$ is also a Wright convex function. However, the converse is not necessarily true always. One of the highly celebrated result was established by C. T. Ng (\cite{Ng}) which states that "A Wright convex function can be decomposed as the algebraic summation of a convex and an additive function. Since then researchers have dealt with different concepts such as characterization, decomposition, approximation for  higher order Wright convex, and approximate-Wright convex functions. More details can can be found in the paper of Maksa and P\'ales \cite{Maksa} and in the book of Kuczma, \cite{Kuczma}. Also in the references therein.\\

Motivated by the definition of Wright convex function, we may think to formulate the notion of Wright convexity in sequence. For all $p,q, r,s\in\N$ with $p<q\leq r<s\in\N$; a sequence $\big<u_n\big>_{n=1}^{\infty}$ is said to be Wright convex if it satisfies the following inequality with alongside mentioned condition
\Eq{*}{
u_q+u_r\leq u_p+u_s \qquad \mbox{provided}\qquad q+r=p+s.
}
Through a characterization, we show that in the case of sequence, a Wright convex and an ordinary convex sequence hold the same meaning.
\section{Main Results}
Throughout this section, $I$ denotes any non-empty  and non-singleton interval in $\R$. For a given $\varepsilon>0$; a function $f:I\to\R$ is said to be $\varepsilon$-affine if for any $x,y\in I$ and for all $t\in[0,1]$; it satisfies the following functional inequality
\Eq{01}{
|f(tx+(1-t)y-tf(x)-(1-t)f(y)|\leq \varepsilon.
}
One can easily observe that if $\varepsilon=0$; then it is nothing but the ordinary affine function.
Before proceeding to the sequence related findings, we first need to recall some results, notions, and terminologies from convex functions and convex geometry.\\

For the function $f:I\to\R$; epigraph is the set defined as
$\{(x,y):f(x)\leq y, x\in I \}\subset{\R^{2}}$. In short we can express this set as $epi(f)$. On the other hand $\R^{2}\setminus epi(f)$ is called hypograph of $f$. If the function $f$ possesses convexity then $epi(f)$ is a convex set. \\

If the intersection of interiors of two convex bodies $A$ and $B(\in\R^{2})$ results in a empty set; then there exists a hyperplane $H(x)=ax+b$ that divides the plane into two half spaces and each of the sets $A$ and $B$ lies in opposite half spaces. This result is known as the separation theorem in convex geometry.\\

\Thm{1}{Let $\varepsilon>0$ and $f:I\to\R$ be a $\varepsilon$-affine function. Then there exists an ordinary affine function $A:I\to\R$ that approximates $f$ such that $||f-A||\leq {\varepsilon}$ holds.}
\begin{proof}
Since the function $f:I\to\R$ is $\varepsilon$-affine; it  satisfies \eq{01}. In other words the following two inequalities together yield \eq{01}
\Eq{*}{
f(tx+(1-t)y)\leq &tf(x)+(1-t)f(y)+{\varepsilon}\\
&\mbox{and}\\
(-f)(tx+(1-t)y)\leq t(-f)(x)&+(1-t)(-f)f(y)+{\varepsilon}
} 
which turn $f$ simultaneously an approximately convex and approximately concave function.

Then as mentioned above; by Hyre's stability theorem there exists a convex $\overline{\C}:I\to\R$ and a concave function $\underline{\C}:I\to\R$ that satisfy the below mentioned inequalities
\Eq{11}{
\overline{\C}(x)-\dfrac{\varepsilon}{2}\leq & f(x)\leq \overline{\C}(x)+\dfrac{\varepsilon}{2}\\
&\mbox{and}\\
\underline{\C}(x)-\dfrac{\varepsilon}{2}\leq &f(x)\leq \underline{\C}(x)+\dfrac{\varepsilon}{2}.\\
}
Combining the above two inequalities, we arrive at
\Eq{*}{
\underline{\C}(x)-\dfrac{\varepsilon}{2}\leq & f(x)\leq \overline{\C}(x)+\dfrac{\varepsilon}{2}}
which shows that the concave function $\underline{\C}-\dfrac{\varepsilon}{2}$ and the convex function 
$\overline{\C}+\dfrac{\varepsilon}{2}$ are not overlapping. This also indicates that the interiors intersection of convex sets $hypo\big(\underline{\C}-\dfrac{\varepsilon}{2}\big)$ and $epi\big(\overline{\C}+\dfrac{\varepsilon}{2}\big)$ is empty.
Then by the separation theorem, there exists a hyperplane $H:\R\to\R$ that separates both of these convex sets. Considering $A:I\to\R$ as the restricted $H\bigg|_{I\cap\R}$; we conclude the following inequality 
\Eq{12}{
\underline{\C}(x)-\dfrac{\varepsilon}{2}\leq  A(x)\leq \overline{\C}(x)+\dfrac{\varepsilon}{2}\quad\quad(x\in I).}
Also, utilizing the first part of the first inequality of \eq{11} together with the last part of the second inequality of \eq{11} we arrive at
\Eq{13}{-\underline{\C}(x)-\dfrac{\varepsilon}{2}\leq - f(x)\leq -\overline{\C}(x)+\dfrac{\varepsilon}{2}\quad\quad(x\in I).}
Summing up \eq{13} together with \eq{12} side by side, we arrive at
\Eq{*}{
-{\varepsilon}\leq A(x)-f(x)\leq {\varepsilon}\quad\quad(x\in I).}
This completes the proof of the statement.
\end{proof}
Before proceeding to the next theorem, we have to go through a fractional inequality first. The detailed proof of the inequality is also included in our recently accepted paper \cite{}.
\Lem{22}{Let $n\in\N$ be arbitrary. Then for any  $a_1,\cdots,a_n\in\R$ and $b_1,\cdots,b_n\in\R_+$; the following inequalities hold
\Eq{111}{
\min\bigg\{\dfrac{a_1}{b_1},\cdots, \dfrac{a_n}{b_n}\bigg\}
\leq\dfrac{a_1+\cdots+a_n}{b_1+\cdots+b_n}\leq \max\bigg\{\dfrac{a_1}{b_1},\cdots, \dfrac{a_n}{b_n}\bigg\}.}}
\begin{proof}
One can easily establish the above inequality by using mathematical induction. Hence, the proof of it is left to the reader.
\end{proof}
In the introduction of the paper, we mentioned about $\varepsilon$-convexity. In the upcoming theorem an equivalent definition of $\varepsilon$-convex function is needed. A function $f:\R_+\to\R$ is said to be $\varepsilon$-convex, if for any $x,u,y\in\R_+$ with $x<u<y$, the following functional inequality holds
\Eq{122}{
\dfrac{f(u)-f(x)-\varepsilon}{u-x}\leq \dfrac{f(y)-f(u)+\varepsilon}{y-u}.
}
\Thm{10}{Let $\big<u_n\big>_{n=0}^{\infty}$ be a $\varepsilon$-convex sequence. Then there exists a ${\varepsilon}$-convex function 
\newline
$f:\R_+\to\R$ such that it satisfies $f(n)=v_n$ for all $n\in\N\cup\{0\}$. Additionally, there exists a convex sequence $\big<v_n\big>_{n=0}^{\infty}$ such that $|u_n-v_n|\leq \varepsilon/2$ holds for all $n\in\N\cup\{0\}$}.
\begin{proof}
To prove the statement, we  define a function $f:\R_+\to\R$ as follows
\Eq{7878}{f(x):=tu_{n-1}+(1-t)u_{n}\quad \mbox{where} \quad x:=t(n-1)+(1-t)n \quad (t\in [0,1]\, \mbox{  and  }\, n\in \N).}
By construction of $f$, it is evident that $f(n)=v_n$ for all $n\in\N$.
Now, We need to prove that $f$ is a $\varepsilon$-convex function. But before that, we need to show two important inequalities. Let $x,y\in\R_+$ with $x<y$ such that $x\in[n_1,n_1+1]$ and $y\in[n_2-1,n_2]$ ($n_1<n_2$ must hold). Then
\Eq{991}{
\qquad\dfrac{f(y)-f(x)-\varepsilon}{y-x}&\leq \max_{n_1\,\leq \,n \,\leq n_2-1}\bigg\{f(n+1)-f(n)-\dfrac{\varepsilon}{n_2-n_1} \bigg\}\\
\newline\\
&\qquad\qquad \mbox{and}\qquad \qquad\qquad\qquad\quad\quad\quad\quad\quad\quad\quad\quad\quad\quad(n\in\N)\\
\newline
\qquad \qquad \min_{n_1\,\leq \,n \,\leq n_2-1}\bigg\{&f(n+1)-f(n)+\dfrac{\varepsilon}{y-x} \bigg\}\leq\dfrac{f(y)-f(x)+\varepsilon}{y-x}.
}
With our assumptions regarding $x$ and $y$, we can rewrite the expression $\dfrac{f(y)-f(x)-\varepsilon}{y-x}$ as follows
\Eq{*}{
\dfrac{\bigg(f(y)-f(n_2-1)\bigg)+\cdots+\bigg(f(n_1+1)-f(x)\bigg)}{\big(y-(n_2-1)\big)+\cdots+\big(n_1+1-x)}-\dfrac{\varepsilon}{y-x}.
}
Using \eq{111} of \lem{22} and utilizing the basic slope property of a line segment, we can conclude 

\Eq{*}{
\dfrac{f(y)-f(x)-\varepsilon}{y-x}&\leq \max\left\{\dfrac{f(y)-f(n_2-1)}{y-(n_2-1)},\cdots \cdots,\dfrac{f(n_1+1)-f(x)}{(n_1+1)-x}
\right\}-\dfrac{\varepsilon}{y-x}\\
&\leq \max\left\{\dfrac{f(n_2)-f(n_2-1)}{n_2-(n_2-1)},\cdots \cdots,\dfrac{f(n_1+1)-f(n_1)}{(n_1+1)-n_1}
\right\}-\dfrac{\varepsilon}{n_2-n_1}\\
&=\max\bigg\{f(n_2)-f(n_2-1)-\dfrac{\varepsilon}{n_2-n_1}, \cdots\cdots, f(n_1+1)-f(n_1)-\dfrac{\varepsilon}{n_2-n_1}\bigg\}.}
This establishes the first part of \eq{991}.
Similarly to show the second inequality of \eq{991}, we proceed by utilizing \eq{111} of \lem{22} as follows 

\Eq{*}{\dfrac{f(y)-f(x)+\varepsilon}{y-x}&\geq \min\left\{\dfrac{f(y)-f(n_2-1)}{y-(n_2-1)},\cdots \cdots,\dfrac{f(n_1+1)-f(x)}{(n_1+1)-x}
\right\}+\dfrac{\varepsilon}{y-x}\\
&=\min\left\{\dfrac{f(n_2)-f(n_2-1)}{n_2-(n_2-1)},\cdots \cdots,\dfrac{f(n_1+1)-f(n_1)}{(n_1+1)-n_1}\right\}+\dfrac{\varepsilon}{y-x}\\
&=\min\left\{f(n_2)-f(n_2-1)+\dfrac{\varepsilon}{y-x},\cdots, f(n_1+1)-f(n_1)+\dfrac{\varepsilon}{y-x}\right\}.}
This validates the second inequality of \eq{991} and completes the establishment. \\
  
Now to justify $\varepsilon$-convexity of $f$, we assume that $x,$ $u$, and $y\in\R_+$ with $x<u<y$ along with the conditions $x\in[n_1,n_1+1]$  $u\in[n_2-1,n_2]$ and $y\in[n_3-2,n_3-1]$ where $n_1,n_2,n_3\in\N$. This also indicates $n_1+1\leq n_2\leq n_3-1$ and yields that $x,$ $u$ and $z$ can lie in all the possible combinations of intervals without affecting the order.  Then according to the results obtained in\eq{991}; we will have atleast one $i\in[n_1,n_2[\cap\N$ and atleast one $j\in [n_2,n_3-1[\cap\N$ such that the two following inequalities are satisfied
\Eq{990}{
\dfrac{f(u)-f(x)-\varepsilon}{u-x}&\leq f(i+1)-f(i)-\dfrac{\varepsilon}{n_2-n_1}\\
&\mbox{and}\\
f(j+1)-f(j)+\dfrac{\varepsilon}{y-u}&\leq \dfrac{f(y)-f(u)+\varepsilon}{y-u}.
}
Due to $\varepsilon$-convexity of the sequence $\Big<u_n\Big>_{n=0}^{\infty}$ $\bigg($or
$\Big<f(n)\Big>_{n=1}^{\infty}\bigg)$; the inequality \eq{1} can be extended to the following extended form 
\Eq{*}{
f(i+1)-f(i)-\dfrac{\varepsilon}{n_2-i}&\leq f(i+1)-f(i)-\dfrac{\varepsilon}{n_2-n_1}\\
&\leq f(j+1)-f(j)\leq f(j+1)-f(j)+\dfrac{\varepsilon}{y-u}.}
With the help of the above inequality, inequalities in \eq{990} can be combined and elaborated as
\Eq{*}{
\dfrac{f(u)-f(x)-\varepsilon}{u-x}&\leq f(i+1)-f(i)-\dfrac{\varepsilon}{n_2-n_1}\leq
f(j+1)-f(j)+\dfrac{\varepsilon}{y-u}&\leq \dfrac{f(y)-f(u)+\varepsilon}{y-u}.
}
This establishes \eq{122}.
Since $x,u$, and $y\in\R_+$ are arbitrary; it yields that the function $f$ possesses $\varepsilon$-convexity. This completes the first part of the theorem.

To prove the second part, we will utilize one of the decomposition results of functional $\varepsilon$-convexity that was formulated by {Hyers} and Ulam in their paper \cite {Hyers}. In our context, it can be formulated as follows

"If a function $f:\R_+\to\R$ is $\varepsilon$-convex then there exists a convex function $g:\R_+\to\R$ that satisfies the inequality
\Eq{*}
{|f(x)-g(x)|\leq \dfrac{\varepsilon}{2}\quad\mbox{for all}\quad x\in I."}
Now we choose the elements of the sequence $\big<v_n\big>_{n=0}^{\infty}$ as $v_n=g(n)$ for all $n\in\N\cup\{0\}$. This yields that $\big<v_n\big>_{n=0}^{\infty}$ is a convex sequence. Using the stated result above, we can obtain the inequality $|f(n)-g(n)|\leq \dfrac{\varepsilon}{2}$; equivalently $|u_n-v_n|\leq \dfrac{\varepsilon}{2}$ for all $n\in\N\cup\{0\}$. It completes the proof.
\end{proof}
The theorem below is a direct consequence of \thm{1} and \thm{10}. Hence, we only provide a scratch of the proof.
\Thm{11}{Let $\big<u_n\big>_{n=0}^{\infty}$ be a $\varepsilon$-affine sequence. Then there exists an arithmetic sequence  $\big<a_n\big>_{n=0}^{\infty}$ such that  $|u_n-a_n|\leq \varepsilon$ holds for all $n\in\N\cup\{0\}$.}
\begin{proof}
We can proceed by constructing $f$ as in \eq{7878}. Then we establish that both $f$ and $-f$ are $\varepsilon$-convex. Utilizing \thm{1}; we are assured that there exists a affine function $A:\R_+\to\R$ such that the inequality $||f-A||\leq\varepsilon$ holds. Extracting the values of these two functions at the points $\N\cup\{0\}$, we have $|f(n)-A(n)|\leq\varepsilon$. Clearly, $\Big<f(n)\Big>_{n=0}^{\infty}$ denotes the sequence $\big<u_n\big>_{n=0}^{\infty}$; while $\big<A(n)\big>_{n=0}^{\infty}$ stands for an arithmetic progression $\big<a_n\big>_{n=0}^{\infty}$. Thus from the already proved norm bound result of \thm{1}, we can conclude $||u_n-a_n||\leq \varepsilon$ for all $n\in\N\cup\{0\}$.
\end{proof}
In the next theorem, we provide a characterization of convex sequences which shows that unlike function  theory ordinary sequential convexity and Wright convexity don't hold different meanings.
\Thm{09}{
A sequence $\big<u_n\big>_{n=1}^{\infty}$ is convex if and only if for any $p,q,r,s\in\N$ with $p<q\leq r<s$, it satisfies the following discrete functional inequality
\Eq{1111}{
u_q+u_r\leq u_p+u_s \qquad \mbox{provided}\qquad q+r=p+s.
}
}
\begin{proof}
At first, we assume that the sequence $\big<u_n\big>_{n=1}^{\infty}$ is convex. This leads us to the following system of inequalities
\begin{equation}
\tag{P}
{ u_{p+1}-u_p\leq u_{p+2}-u_{p+1}}
\end{equation}
\begin{equation}
\tag{P+1}
{ u_{p+2}-u_{p+1}\leq u_{p+3}-u_{p+2}}
\end{equation}
$$\vdots \qquad \qquad \vdots$$
\begin{equation}
\tag{S-3}
{ u_{s-2}-u_{s-3}\leq u_{s-1}-u_{s-2}}
\end{equation}
\begin{equation}
\tag{S-2}
{ u_{s-1}-u_{s-2}\leq u_{s}-u_{s-1}}.
\end{equation}
From this system of inequalities, summing up all side by side, we obtain 
\Eq{7788}{
u_{p+1}+u_{s-1}\leq u_p+u_s.
} 
Again adding up all inequalities side by side excluding the inequalities $P$ and $S-2$, we arrive at
\Eq{7789}{
u_{p+2}+u_{s-2}\leq u_{p+1}+u_{s-1}.
}
Similarly  proceeding as above, excluding the inequalities $P$, $P+1$, $S-3$, and $S-2$, we add all remaining inequalities to get
\Eq{7790}{
u_{p+3}+u_{s-3}\leq u_{p+2}+u_{s-2}.
}
And we keep continuing the process and each time we find new inequalities similar to \eq{7788}, \eq{7789} and \eq{7790}. Depending upon the values of $p$ and $q$, we will end up in a system of inequalities that can be summarized as follows
\Eq{*}{
2u_{\frac{p+s}{2}}\leq u_{\frac{p+s}{2}-1}+ u_{\frac{p+s}{2}+1}\leq &\cdots \leq u_{p+1}+u_{s-1}\leq u_{p}+u_{s} \quad\mbox{if $p+s$ is even}\\
&\mbox{or}\\
u_{\frac{p+s-1}{2}}+ u_{\frac{p+s+1}{2}}\leq u_{\frac{p+s-1}{2}-1}+ u_{\frac{p+s-1}{2}+1}\leq&\cdots \leq u_{p+1}+u_{s-1}\leq u_{p}+u_{s} \quad\mbox{if $p+s$ is odd.}
}
The above inequalities show all the 
possible combinations of $q,r\in\N$ such that  \eq{1111} is satisfied.
This establishes the first part of our theorem.\\
 
Conversely, we assume that for the sequence $\big<u_n\big>_{n=0}^{\infty}$, inequality \eq{1111} holds. Let $n\in\N\setminus\{1\}$ be arbitrary. Substituting 
$p=n-1,$ $ q=r=n$ and $s=n+1$ in \eq{1111}; we arrive at \eq{99}. This shows that $\big<u_n\big>_{n=1}^{\infty}$ possesses sequential convexity and completes the proof.
\end{proof}
In this paper, we dealt with only one specific type of approximately convex sequence. One can investigate another types of such sequences and find out the proper decomposing formulation and interesting results.
\bibliographystyle{plain}

\end{document}